# Planning Skip-Stop Transit Service under Heterogeneous Demands

Yu Mei [a], Weihua Gu [a], Michael Cassidy [b], Wenbo Fan [c*]

[a] Department of Electrical Engineering, The Hong Kong Polytechnic University, Hung Hom, Kowloon, Hong Kong SAR, China
[b] Department of Civil and Environmental Engineering, University of California, Berkeley, USA
[c] School of Transportation and Logistics, Southwest Jiaotong University, Chengdu, China

**ABSTRACT**

Transit vehicles operating under skip-stop service visit only a subset of the stops residing along a corridor. It is a strategy commonly used to increase vehicle speeds and reduce patron travel times. The present paper develops a continuous approximation model to optimally design a select form of skip-stop service, termed AB-type service. The model accounts for spatially-heterogeneous demand patterns. An efficient heuristic is developed to obtain solutions. These are shown to be near-optimal for a variety of numerical examples. Results also indicate that optimal AB-type designs outperform optimized all-stop service in a variety of cases. The AB-type service is found to be especially competitive when travel demands are high, trip origins are unevenly distributed along a corridor, and patrons have relatively high values of time. In these cases, AB-type service is found to reduce system costs by as much as 8%.

**Keywords**: transit corridor design; skip-stop service; continuous approximation; heterogeneous demand

[*] Corresponding author.
E-mail address: wbfan@swjtu.edu.cn.

# 1 Introduction

Skip-stop service can be used when multiple transit lines operate simultaneously in a single corridor. With this service, the buses or trains assigned to each line visit only a subset of the corridor's stops (Vuchic, 2005; Leiva et al., 2010; Chen et al., 2015b; Larrain and Muñoz, 2016; Soto et al., 2017). By skipping certain stops, each transit vehicle can operate at higher commercial speeds than would be achieved under traditional, all-stop service. Skip-stop service can also save patrons' travel times, though this benefit diminishes when patrons are forced to transfer between distinct lines.

Our present interest lies in a particular form of skip-stop service, commonly known as AB-type service. The name dates back to the middle of the last century, when Chicago's metro-rail system operated two parallel lines (designated A and B) in select corridors. This same service form was implemented in the Metro system in Santiago, Chile (Metro de Santiago, 2008), and in bus corridors in Portland, Oregon (Trimet.org, n.d.). The idea is illustrated in Figure 1. Vehicles serving lines A and B visit non-transfer stops (shown as unshaded squares in the figure) in alternating fashion, as well as every transfer stop (darkened circles). A patron who accesses this service via a Line-A (non-transfer) stop and egresses at a Line-B stop must transfer somewhere along her journey.

To clarify by example, two possible means of transferring are shown with dotted arrows in Figure 1. Each is an option when the patron's origin and destination stops belong to two distinct lines, and where both the origin and destination stops reside between the same two transfer stops. (The line segment spanning two transfer stops will be termed a "skip-stop bay.") For route option 1, the patron alights at the first downstream transfer stop, then backtracks to her destination. For route option 2, she backtracks first. The patron will choose the lower-cost route between the two options. Transfers can be achieved without backtracking when origin and destination stops reside in distinct bays.

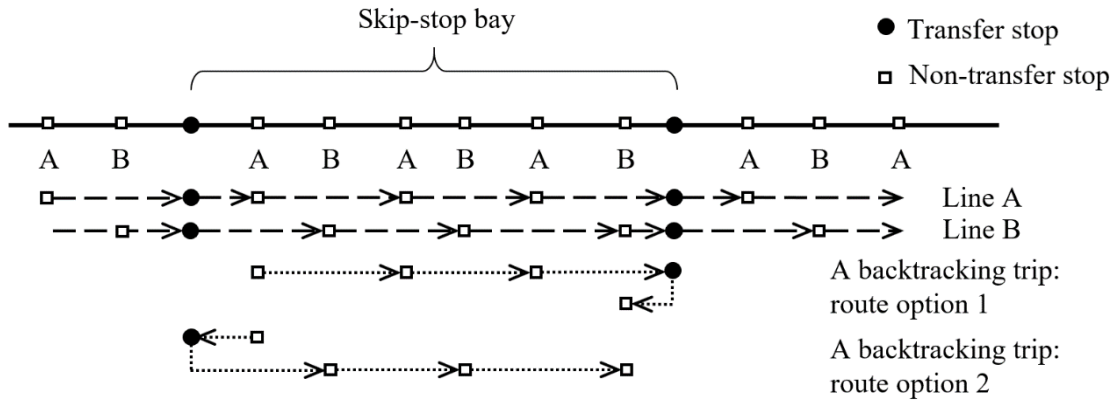

**Figure 1. A typical AB-type service**

This strategy enables vehicles serving distinct lines to travel in a common direction without need to overtake other transit vehicles. An example is shown in Figure 2 for a case of AB-type service with 3 distinct lines. Note how the need for overtaking is averted by dispatching vehicles in such way that the headway is never less than a minimum threshold, $H_{min}$.[1] This means that all lines in a single direction can be economically served via a single bus lane or rail track.

In several previous studies of AB-type service, case-specific models were developed in which inputs and outputs were defined in discrete and detailed fashion (Suh et al., 2002; Lee et al., 2014; Abdelhafiez et al., 2017). Each such model typically generates an optimal number of lines and the set of stops visited by each line in a corridor, where all stops in the corridor are specified a priori. (The number of lines and the set of stops visited by each line will henceforth be termed the "route plan.") To their credit, these models can account for spatially-heterogeneous demands. They cannot, however, be used to optimize a corridor's number of stops. The complexity of these models, moreover, has typically necessitated reliance on heuristically-generated solutions, with optimality gaps that can be difficult to evaluate (Lee et al., 2014; Abdelhafiez et al., 2017). Of perhaps even greater concern, the case-specific nature of these models tends to produce few general insights into relations between model inputs and outputs.

---

[1] The $H_{min}$ might be established with safety considerations in mind; see Gu et al. (2016).

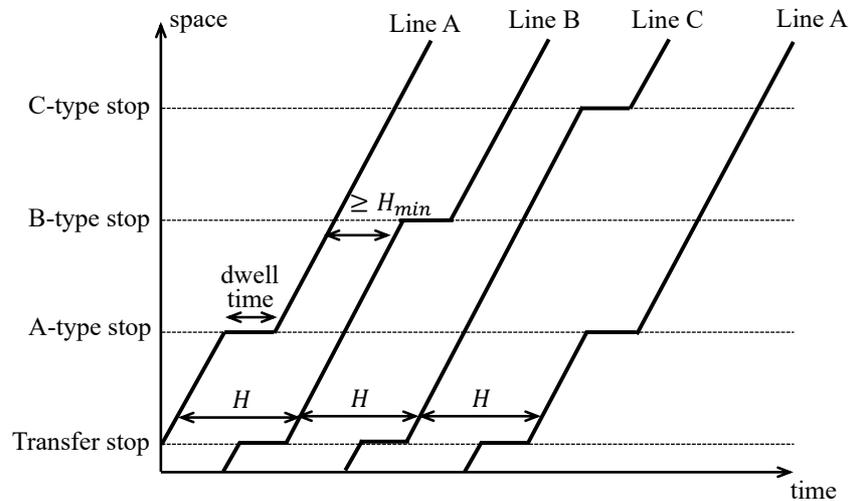

**Figure 2. Transit vehicle trajectories of a three-line AB-type system (Gu et al., 2016)**

Efforts to remedy these limitations have focused on models in which inputs and outputs are specified in continuous, less-detailed fashion. This would seem a promising approach for designing AB-type systems, since continuous and continuous approximation (CA) models can produce both, general insights into cause and effect, and designs that are well-suited to high-level transit planning. The designs are invariably developed for idealized settings, but can be adjusted to suit the peculiarities of real-world environments; see Estrada et al. (2011).

Of particular relevance are the works of Freyss et al. (2013) and Gu et al. (2016). Both developed models that jointly optimize the spacing between stops (and thus the number of stops in a corridor) and the route plan for AB-type systems. Both models unveil insights regarding cause and effect. Yet, both assume that travel demand is distributed uniformly along a corridor. This is unlikely to hold in real settings, and assumes-away conditions for which an AB-type strategy can be especially well suited.

The literature holds clues on how to resolve this limitation as well. There we find several CA models for designing various forms of transit service (though not for AB-type service!) in the presence of spatially-heterogeneous demands (Vaughan and Cousins, 1977; Wirasinghe and Ghoneim, 1981; Vaughan, 1986; Medina et al., 2013; Ouyang et al., 2014). Each model uses multiple continuous functions. Each approximates temporally- or spatially-varying inputs to transit design, such as demand densities; and decision variables (i.e. outputs), such as line and stop spacings. The models were typically solved by decomposing the problems into sub-problems, with each holding only a handful of decision variables. Global optima and even closed-form solutions could often be developed as a result (e.g., Wirasinghe and Ghoneim, 1981). In another relevant work (Chen et al., 2015a), a more complex CA model for designing transit networks with hybrid structures was solved approximately via heuristics with efficient search methods. Optimality gaps were evaluated in that work by comparing solutions against lower bounds.

The present paper continues in this tradition. It has developed a CA model to jointly optimize stop locations and route plans for AB-type transit service in corridors with spatially-heterogeneous demands. With this new model, the stop spacing and the number of stops within a skip-stop bay are both allowed to vary along the corridor. Along shorter segments, lines in a common direction roughly share the same density of stops.

Models that address the challenges of spatially-heterogeneous demands are formulated in the following section. A solution approach that uses calculus of variations is developed in Section 3, and a method to assess optimality gaps in Section 4. Numerical analyses are presented in Section 5, and practical implications in Section 6. Notations used in these efforts are defined when they are introduced, and are also tabulated in Appendix A for the reader's convenience.

## 2 Models

The framework or set-up for our modeling efforts is described in Section 2.1. Means to quantitatively express heterogeneous demand are presented in Section 2.2. Models of user and agency costs are presented in Sections 2.3 and 2.4, respectively. Means of optimizing service designs are furnished in Section 2.5.

### 2.1 Framework

We consider bi-directional service along a closed-loop corridor of length *L*, as exemplified in Figure 3a.[2]

Four scalar variables are used in designing AB-type service. These are: the numbers of skip-stop lines that serve travel demand in the clockwise (CW) and counterclockwise (CCW) directions, denoted $m_{cw}$ and $m_{ccw}$, respectively; and the transit vehicles' directional headways, $H_{cw}$ and $H_{ccw}$. Vehicles traveling CW on each line are dispatched at a common headway, $m_{cw}H_{cw}$, because vehicles serving distinct lines are consecutively dispatched at $H_{cw}$. Vehicles traveling CCW on the same line are dispatched at $m_{ccw}H_{ccw}$.

Two continuous, slowly-varying functions are also needed for our descriptions. The first of these, denoted $T(x)$, approximates the integer-valued number of stops in a skip-stop bay at location $x$ for $0 \le x < L$. The second, denoted $s(x)$, approximates a step-wise function that defines the location-dependent spacing between stops.

All-stop service can be viewed as a special case of the AB-type in which $m_{cw} = m_{ccw} = 1$, and $T(x) = 1, \forall x$. Vehicles traveling in the CW and CCW directions are assumed to visit the same set of stops along the loop, just as in Gu et al. (2016).

Four additional commonly-adopted assumptions are used in the present work. These are: (i) a patron accesses and egresses the system via stops nearest her origin and destination; (ii) a transit vehicle loses a fixed time, $\tau$, at each stop owing to boarding and alighting patrons and to acceleration and deceleration; (iii) patrons arrive randomly to their origin stops, without reference to service schedules, as occurs when vehicle headways are small; and (iv) transit vehicles maintain regular headways.

### 2.2 Demand functions and variables

Demand for transit travel along the corridor is assumed invariant to time.[3] Its density (in units of trips/km²/h) is represented by a slow-varying, integrable function, $\lambda(x, y)$, where $x$ and $y$ are the location coordinates of origin and destination stops, respectively. Without loss of generality, $x$ and $y$ are measured in the CW direction along the corridor: $0 \le x, y < L$. We assume a patron always chooses the travel direction that minimizes her travel distance; i.e., a patron travels CW if $0 < y - x \leq \frac{L}{2}$ or $y - x \leq -\frac{L}{2}$, and CCW if $-L/2 < y - x \leq 0$ or $y - x > L/2$.[4]

---

[2] A loop was chosen to avoid certain operational constraints that are often required of linear corridors, e.g., that the headways for the two service directions must be equal. A loop corridor also allows trip length distributions to be set independently of trip origins. This enables separate parametric analyses of trip length distribution and spatial distribution of trip origins. Our models can readily be adapted to suit linear corridors; see Footnote 5.
[3] The proposed model can be made to account for time-varying demand over a day in a manner similar to Daganzo (2010).

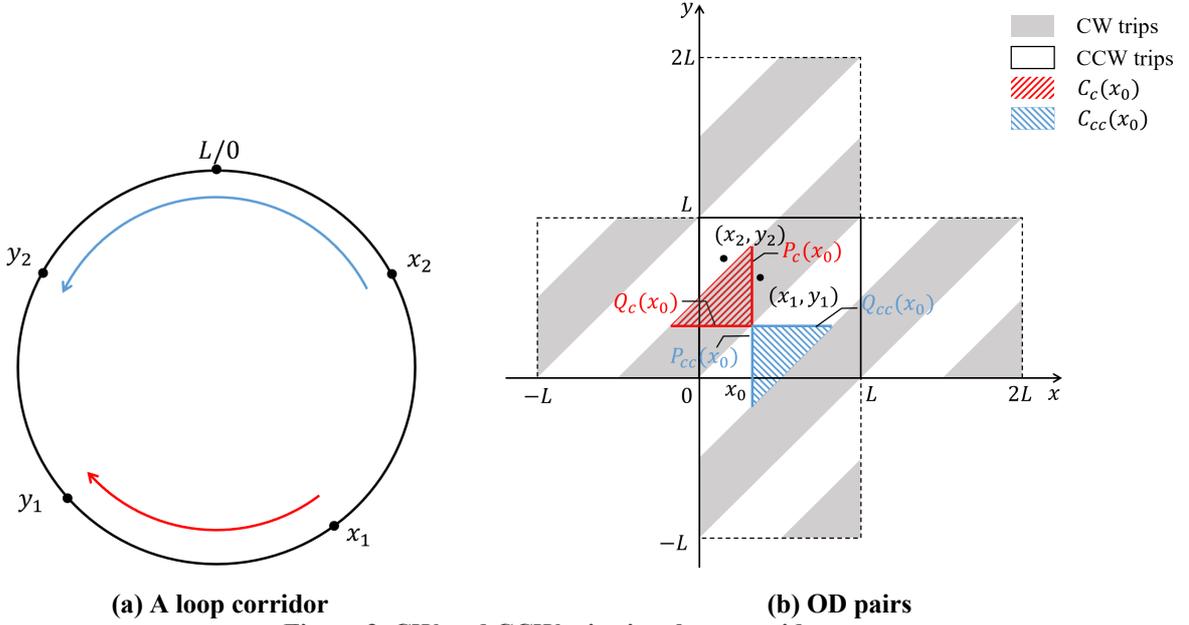

**(a) A loop corridor**          **(b) OD pairs**
**Figure 3. CW and CCW trips in a loop corridor**

A square of size $[0, L] \times [0, L]$ represents the set of OD pairs in the corridor, as shown in the center of Figure 3b. The set of ODs to be served by CW trips are highlighted by the figure's shaded diagonal swaths. The OD pairs served by CCW trips are unshaded. Note as examples the CW trip $(x_1, y_1)$ shown in both Figures 3a and b, and the CCW trip $(x_2, y_2)$ shown in both figures as well.

We define below three aggregate demand functions and variables labeled (i)-(iii) to be used in the present models. To simplify things, the domain of demand function $\lambda(\cdot)$ is expanded to allow $\lambda(x \pm L, y) \equiv \lambda(x, y \pm L) \equiv \lambda(x, y)$ for $0 \leq x, y < L$. Note how this expansion is displayed by the additional squares etched with dashed lines in Figure 3b.

(i) The densities of origins for CW and CCW trips at location $x$ are denoted $P_{cw}(x)$ and $P_{ccw}(x)$, respectively; and the densities of destinations at location $y$ are $Q_{cw}(y)$ and $Q_{ccw}(y)$. We thus have

$$P_{cw}(x) = \int_{y=x}^{x+L/2} \lambda(x,y) dy, P_{ccw}(x) = \int_{y=x-L/2}^{x} \lambda(x,y) dy, 0 \leq x < L, \quad (1a)$$

$$Q_{cw}(y) = \int_{x=y-L/2}^{y} \lambda(x,y) dx, Q_{ccw}(y) = \int_{x=y}^{y+L/2} \lambda(x,y) dx, 0 \leq y < L, \quad (1b)$$

with: $\lambda(x, y)$ defined over the $(x, y)$ plane in Figure 3b; $P_{cw}(\cdot)$ and $Q_{cw}(\cdot)$ obtained by integrating $\lambda(x, y)$ along the vertical and horizontal edges of a triangle, like the upper-left hatched one in the figure; and $P_{ccw}(\cdot)$ and $Q_{ccw}(\cdot)$ similarly obtained, as exemplified by the vertical and horizontal edges of the lower-right hatched triangle.

(ii) Total CW and CCW demand, denoted $\Lambda_{cw}$ and $\Lambda_{ccw}$, are given by:

$$\Lambda_{cw} = \int_{x=0}^{L} P_{cw}(x) dx = \int_{y=0}^{L} Q_{cw}(y) dy, \Lambda_{ccw} = \int_{x=0}^{L} P_{ccw}(x) dx = \int_{y=0}^{L} Q_{ccw}(y) dy. \quad (2)$$

---

[4] In reality, some patrons whose trip length is near $\frac{L}{2}$ may choose to travel in the opposite direction, if that would reduce her travel time. That behavior would have only small effect on generalized cost, and is therefore ignored to simplify the modeling work.

The $\Lambda_{cw}$ is obtained by integrating $\lambda(x,y)$ over the shaded area in an $L \times L$ square in Figure 3b, and $\Lambda_{ccw}$ by integrating $\lambda(x,y)$ over the unshaded area of the same square.

(iii) The approximate flows of patrons aboard CW and CCW vehicles at location $x$ are denoted $C_{cw}(x)$ and $C_{ccw}(x)$, respectively. These are defined as:

$$C_{cw}(x) = \int_{z=x-L/2}^{x} \int_{y=x}^{z+L/2} \lambda(z,y) dy dz, C_{ccw}(x) = \int_{z=x}^{x+L/2} \int_{y=z-L/2}^{x} \lambda(z,y) dy dz, 0 \leq x < L. \quad (3)$$

As an example, Figure 3b shows how a $C_{cw}(x_0)$ for an arbitrary location, $x_0$, is obtained by integrating $\lambda(x,y)$ over the hatching of the upper-left triangle. We note that the area above that triangle's horizontal bar in the shaded swath represents trips destined downstream of $x_0$ in the CW travel direction; and the area to the left of its vertical bar in the same swath represents trips originating upstream of $x_0$ in the same direction. The intersection between the above two areas (i.e., the triangle itself) represents trips that traverse $x_0$. The $C_{ccw}(x_0)$ is similarly obtained by integrating over the lower-right hatched triangle in the figure.

Equations in (3) are approximations because flows of on-board patrons are affected by the locations along the route where stops reside. The approximation is a good one when stop spacings are much smaller than average trip length, as is often the case in large cities (Wirasinghe and Ghoneim, 1981).[5]

### 2.3 User costs

The time that users (patrons) collectively spend in the system each hour is dictated by three activities: (i) accessing and egressing stops by walking, denoted $UT_a$; (ii) waiting at origin and transfer stops, $UT_w$; and (iii) traveling aboard vehicles, $UT_v$. A time penalty for the inconvenience of transferring between vehicles, $UT_t$, is also included in our calculations. Formulas for activity items (i)-(iii) and the transfer penalty are presented below.

*2.3.1 Access and egress time*

Since demand is assumed to change slowly along the corridor, a user's average time spent walking along the route to the nearest stop from her origin at $x$ is $\frac{s(x)}{4v_w}$, where $v_w$ is walking speed. The same formula holds for the average egress time from a stop nearest a destination at $x$. The time collectively spent each hour accessing and egressing the system is thus

$$UT_a = \int_{x=0}^{L} \frac{s(x)}{4v_w} \big(P_{cw}(x) + Q_{cw}(x) + P_{ccw}(x) + Q_{ccw}(x)\big) dx. \quad (4)$$

*2.3.2 Wait time*

Five distinct trip types are defined below for CW travel. Similar definitions apply to five CCW trip types, but are omitted for brevity. Wait times vary with trip types, as will become evident.

A type-1 trip is one in which both the origin and destination stops are transfer stops. The average wait time in this case is $W_1 = H_{cw}/2$, since a patron can travel via any CW line without transferring between vehicles. To determine the total number of type-1 trips made per hour, $N_1$, we note that the probability that a trip from $x$ to $y$ ($0 \leq x, y < L$) is of type-1 can be expressed as $\frac{1}{T(x)} \cdot \frac{1}{T(y)}$, since $\frac{1}{T(\cdot)}$ is the probability that the nearest stop to $x$ or $y$ is a transfer stop. Hence, $N_1 = \iint_{D_{cw}} \frac{\lambda(x,y)}{T(x)T(y)} dxdy$, where $D_{cw}$ is the set of CW OD pairs: $D_{cw} \equiv \{(x,y) | x < y \leq x + L/2, x \in [0, L)\}$.

---

[5] If a linear corridor is modeled, one needs only to replace equations (1a), (1b) and (3) respectively by: $P_{cw}(x) = \int_{y=x}^{L} \lambda(x,y) dy$, $P_{ccw}(x) = \int_{y=0}^{x} \lambda(x,y) dy$; $Q_{cw}(y) = \int_{x=0}^{y} \lambda(x,y) dx$, $Q_{ccw}(y) = \int_{x=y}^{L} \lambda(x,y) dx$; and $C_{cw}(x) = \int_{z=0}^{x} \int_{y=x}^{L} \lambda(z,y) dy dz$, $C_{ccw}(x) = \int_{z=x}^{L} \int_{y=0}^{x} \lambda(z,y) dy dz$.

A type-2 trip occurs between a transfer stop and a non-transfer one. The average wait is $W_2 = m_{cw}H_{cw}/2$, since the patron must travel on a specified line. The probability that a trip from $x$ to $y$ is of type-2 is $\frac{1}{T(x)} \cdot \frac{T(y)-1}{T(y)} + \frac{T(x)-1}{T(x)} \cdot \frac{1}{T(y)}$, such that the hourly number of these trips is $N_2 = \iint_{D_{cw}} \frac{(T(x)+T(y)-2)\lambda(x,y)}{T(x)T(y)} dxdy$.

A type-3 trip occurs between two non-transfer stops along a single line. The wait is again $W_3 = m_{cw}H_{cw}/2$. The trip's probability is $\frac{1}{m_{cw}} \frac{T(x)-1}{T(x)} \cdot \frac{T(y)-1}{T(y)}$; and its number is $N_3 = \iint_{D_{cw}} \frac{(T(x)-1)(T(y)-1)\lambda(x,y)}{m_{cw}T(x)T(y)} dxdy$.

A type-4 trip occurs between two non-transfer stops along two distinct lines, when both stops reside in the same skip-stop bay. The trip requires back-tracking, as previously shown with the dotted arrows in Figure 1. Given this requirement, the average wait at origin and transfer stops is $W_4 = (m_{cw}H_{cw} + m_{ccw}H_{ccw})/2$. Since type-4 trips are short in length, their density at $x$, $b_{cw}(x)$, can be approximated as:

$$b_{cw}(x) \approx \frac{(m_{cw}-1)(T(x)-1)^2}{m_{cw}T(x)^3 s(x)} \int_{z=x-\frac{T(x)s(x)}{2}}^{x+\frac{T(x)s(x)}{2}} \int_{y=z}^{x+\frac{T(x)s(x)}{2}} \lambda(z,y) dz dy, \tag{5}$$

where we specify that $\lambda(x + L, y + L) \equiv \lambda(x - L, y - L) \equiv \lambda(x, y)$ for $0 \leq x, y < L$ to ensure that the integral in (5) is well defined. The derivation of (5) is relegated to Appendix B.1. The hourly number of these trips is $N_4 = \int_{x=0}^{L} b_{cw}(x) dx$.

Finally, a type-5 trip occurs between two non-transfer stops along two distinct lines when at least one transfer stop resides between. Since the trip's needed transfer does not entail backtracking, the average wait is $W_5 = m_{cw}H_{cw}$. The hourly number of trips is $N_5 = \Lambda_{cw} - \sum_{i=1}^{4} N_i$.

The above results are summarized in Appendix B.2 for readers' ease of reference.

The wait collectively incurred each hour among all CW trips is $\sum_{i=1}^{5} N_i W_i$. By combining this expression with that for CCW trips, one obtains

$$UT_w = \frac{(2m_{cw}-1)H_{cw}\Lambda_{cw}}{2} + \frac{(2m_{ccw}-1)H_{ccw}\Lambda_{ccw}}{2} + \int_{x=0}^{L} \left[ -\frac{(m_{cw}-1)H_{cw}}{2} \frac{P_{cw}(x)+Q_{cw}(x)}{T(x)} - \frac{(m_{ccw}-1)H_{ccw}}{2} \frac{P_{ccw}(x)+Q_{ccw}(x)}{T(x)} + \frac{m_{ccw}H_{ccw}-m_{cw}H_{cw}}{2} \left( b_{cw}(x) - b_{ccw}(x) \right) \right] dx, \tag{6}$$

where $b_{ccw}(x)$ is the approximate density of CCW trips involving backtracking, as obtained by replacing $m_{cw}$ with $m_{ccw}$ in (5). See Appendix B.2 for the derivation of (6).

*2.3.3 In-vehicle travel time*

The commercial speed of a CW vehicle in the vicinity of $x$ is

$$v_{cw}(x) = \left[ v^{-1} + \frac{(T(x)-1)/m_{cw}+1}{T(x)s(x)} \tau \right]^{-1}, \tag{7}$$

where: $v^{-1}$ is the time to travel a unit distance at cruise speed $v$; and $\frac{(T(x)-1)/m_{cw}+1}{T(x)s(x)}$ is the number of stops visited by the vehicle per unit distance. The commercial speed of a CCW vehicle, $v_{ccw}(x)$, is obtained by replacing $m_{cw}$ in (7) with $m_{ccw}$.

In-vehicle travel time for all CW trips combined consists of the time collectively spent traveling in the CW direction, $\int_{x=0}^{L} C_{cw}(x) v_{cw}(x)^{-1} dx$; and backtracking. The latter is obtained knowing the extra distance travelled due to backtracking near $x$. This is roughly approximated as $T(x)s(x)/3$; see Gu et al. (2016).[6] Since half the extra distance is travelled CW and the other half CCW, the average extra time that a user spends per backtracking trip is $\frac{1}{6}T(x)s(x) \cdot (v_{cw}(x)^{-1} + v_{ccw}(x)^{-1})$; and the extra time collectively spent each hour by all users is $\int_{x=0}^{L} \frac{1}{6}T(x)s(x) \cdot (v_{cw}(x)^{-1} + v_{ccw}(x)^{-1}) b_{cw}(x) dx$.

Combining the above with the in-vehicle time for CCW trips yields

$$UT_v = \int_{x=0}^{L} C_{cw}(x)\left(v^{-1} + \frac{(T(x)-1)/m_{cw}+1}{T(x)s(x)}\tau\right) + C_{ccw}(x)\left(v^{-1} + \frac{(T(x)-1)/m_{cw}+1}{T(x)s(x)}\tau\right) + (b_{cw}(x) + b_{ccw}(x))\left[\frac{T(x)s(x)}{3v} + \frac{\tau}{6}\left(\frac{T(x)-1}{m_{cw}} + \frac{T(x)-1}{m_{ccw}} + 2\right)\right] dx. \quad (8)$$

*2.3.4 Transfer penalty*

Only trip types 4 and 5 entail transfers; one per trip. The hourly number of transfers across all CW trips is thus $\iint_{D_{cw}} \lambda(x,y) \frac{m_{cw}-1}{m_{cw}} \frac{T(x)-1}{T(x)} \cdot \frac{T(y)-1}{T(y)} dxdy$. To enable the development of an efficient solution approach via calculus of variations, this double integral is approximated as a single one in Appendix B.3. By combining the CW and CCW directions, the hourly transfer penalty becomes

$$UT_t \approx C_t \left[\frac{m_{cw}-1}{m_{cw}}\Lambda_{cw} + \frac{m_{ccw}-1}{m_{ccw}}\Lambda_{ccw} - \int_{x=0}^{L}\left[\frac{m_{cw}-1}{m_{cw}}(P_{cw}(x) + Q_{cw}(x)) + \frac{m_{ccw}-1}{m_{ccw}}(P_{ccw}(x) + Q_{ccw}(x))\right]\frac{(2T(x)-1)}{2T^2(x)}dx\right], \quad (9)$$

where $C_t$ is the penalty cost per transfer in a unit of time.

## 2.4 Agency cost

The hourly cost incurred by the transit agency consists of four components: (i) a distance-based vehicle operating cost, primarily fuel cost, denoted $AC_K$; (ii) a time-based vehicle operating cost, including the amortized vehicle purchase cost and staff wages, $AC_H$; (iii) an amortized infrastructure cost for bus lanes or rail tracks, $AC_I$; and (iv) an amortized infrastructure cost for stops, $AC_S$. The derivations for these costs are simple and similar to those in Daganzo and Ouyang (2019) and Gu et al. (2016), and are thus omitted for brevity. They are expressed as

$$AC_K = \frac{\pi_k L}{\mu}\left(\frac{1}{H_{cw}} + \frac{1}{H_{ccw}}\right) \quad (10)$$

$$AC_H = \frac{\pi_h}{\mu} \int_{x=0}^{L} \left(\frac{1}{H_{cw}} v_{cw}(x)^{-1} + \frac{1}{H_{ccw}} v_{ccw}(x)^{-1}\right) dx \quad (11)$$

$$AC_I = \frac{2\pi_i L}{\mu} \quad (12)$$

$$AC_S = \frac{\pi_s}{\mu} \int_{x=0}^{L} \frac{1}{s(x)} dx, \quad (13)$$

where: $\pi_k$ and $\pi_h$ are the operating costs per vehicle-km and per vehicle-hour of service, respectively; $\pi_i$ and $\pi_s$ are the amortized hourly construction and maintenance costs per km of transit line and per stop; and $\mu$ is the patrons' value of time, for which the average hourly wage rate might be used as a proxy.

---

[6] This is a coarse approximation despite the assumption of slow-varying demand. However, the error brought by this assumption is small, since the backtracking costs are small under an optimal design.

## 2.5 Joint optimization problem

The problem of minimizing the generalized cost of service, $GC$, can be formulated as

$$\min_{m_{cw},m_{ccw},H_{cw},H_{ccw},s(x),T(x)} GC = UT_a + UT_w + UT_v + UT_t + AC_K + AC_H + AC_I + AC_S \quad (14a)$$

subject to:
$$m_{cw}, m_{ccw} \in \{1,2,3,4\} \quad (14b)$$

$$H_{cw} \geq H_{min} + I(m_{cw} > 1) \cdot \tau \quad (14c)$$

$$H_{ccw} \geq H_{min} + I(m_{ccw} > 1) \cdot \tau \quad (14d)$$

$$H_{cw} \leq \frac{K}{\max_{0 \leq x < L}\{C_{cw}(x) + T(x)s(x)(b_{cw}(x) + b_{ccw}(x))/2\}} \quad (14e)$$

$$H_{ccw} \leq \frac{K}{\max_{0 \leq x < L}\{C_{ccw}(x) + T(x)s(x)(b_{cw}(x) + b_{ccw}(x))/2\}} \quad (14f)$$

$$s(x) > 0, T(x) \geq 1, \forall x \in [0, L), \quad (14g)$$

where: $H_{min}$ was previously defined with the aid of Figure 2; $I(\cdot)$ is an indicator function that returns 1 if the argument is true, and 0 otherwise; and $K$ is a vehicle's passenger-carrying capacity. We note that $m_{cw}$ or $m_{ccw} = 1$ denote all-stop service for the CW or CCW directions.

Constraint (14b) allows for a maximum of four lines in each direction, since a larger number would likely be difficult for users to navigate. Constraints (14c-d) maintain minimum service headways (e.g. for safety), whether service is AB-type or all-stop.

Constraint (14e) limits CW patron flows based on vehicle capacities. The term $T(x)s(x)(b_{cw}(x) + b_{ccw}(x))/2$ is the added on-board flow due to backtracking. We note that the total flow for both CW and CCW trips with backtracking is $T(x)s(x)(b_{cw}(x) + b_{ccw}(x))$, half of which is added to the CW travel direction. Constraint (14f) is similarly derived for CCW flows. Both constraints assume that patrons on type-5 trips randomly choose transfer stops between their origins and destinations, such that the flows of these patrons are evenly distributed across all lines.

Constraints (14g) are boundary constraints for the continuous functions $s(x)$ and $T(x)$.

The complexity of (14a-g) renders exact solutions out of reach. A heuristic approach to solving the CA model is presented below.

## 3 Heuristic Solution Method

The components of (4), (6) and (8)-(13) enable a re-writing of (14a) as follows:

$$GC = h(m_{cw}, m_{ccw}, H_{cw}, H_{ccw}) + \int_{x=0}^{L} G(m_{cw}, m_{ccw}, H_{cw}, H_{ccw}, s(x), T(x), x)dx, \quad (15)$$

where: $h$ is a function solely of scalar variables; and $G$ is a function of both scalar variables and continuous decision functions.

We propose an iterative algorithm to minimize $GC$ in two stages. These are separately described in Sections 3.1 and 3.2.

### 3.1 Stage 1: finding $s^*(x)$ and $T^*(x)$

The first stage entails determining the decision functions $s(x)$ and $T(x)$ to minimize $GC$ for fixed values of the scalar variables in (15), and subject to constraints (14e-g). Calculus of variations is used to render the minimization of $\int_{x=0}^{L} G(\cdot, x)dx$ equivalent to minimizing $G(\cdot, x)$ for each $x \in [0, L]$. This is done by selecting $n$ discrete points, denoted $x_j \equiv (j - 0.5)\Delta x, j = 1,2,...,n$, from the continuous range of $x$ at equal spacing, $\Delta x = L/n$. The values of $s^*(x_j)$ and $T^*(x_j)$ are found to minimize $G(\cdot, x_j)$ for each $x_j$.

The minimization is a complex task, because determining $b_{cw}(x_j)$ and $b_{ccw}(x_j)$ each entails double integration; see again (5). An iterative search method is proposed to tackle this complexity.

In each iteration, the values of $b_{cw}(x_j)$, $b_{ccw}(x_j)$ and $T(x_j)$ are fixed. The $s^*(x_j)$ is determined as:

$$s^*(x_j) = \min\left\{\tilde{s}(x_j), \frac{\min\{2K/H_{cw}-2C_{cw}(x_j), 2K/H_{ccw}-2C_{ccw}(x_j)\}}{(b_{cw}(x_j)+b_{ccw}(x_j))T(x_j)}\right\}, \quad (16a)$$

where: the second argument comes from the vehicle capacity constraints (14e-f); and the first argument, $\tilde{s}(x_j)$, is derived from the first-order condition of $G(\cdot, x)$ with respect to $s(x)$. Thus,

$$\tilde{s}(x) = \sqrt{\frac{\left(\frac{1}{m_{cw}}+\frac{m_{cw}-1}{m_{cw}T(x)}\right)\left(C_{cw}(x)+\frac{\pi_h}{\mu H_{cw}}\right)\tau+\left(\frac{1}{m_{ccw}}+\frac{m_{ccw}-1}{m_{ccw}T(x)}\right)\left(C_{ccw}(x)+\frac{\pi_h}{\mu H_{ccw}}\right)\tau+\frac{\pi_s}{\mu}}{(P_{cw}(x)+Q_{cw}(x)+P_{ccw}(x)+Q_{ccw}(x))/4v_w+T(x)(b_{cw}(x)+b_{ccw}(x))/3v}}. \quad (16b)$$

Equations (16a-b) are valid because $G(\cdot, x)$ is a linear combination of $s(x)$ and $s(x)^{-1}$ with non-negative coefficients. (Readers can verify this by checking each cost component in (4-13) that involves $s(x)$.) Thus, $G(\cdot, x)$ is a convex function of $s(x)$.

We then enumerate $T(x_j) \in \{1,2,...,30\}$ for each $j$ to find the optimal skip-stop bay sizes, $T^*(x_j)$.[7]

Finally, the $b_{cw}(x_j)$ and $b_{ccw}(x_j)$ are updated using the method of successive averages (Sheffi, 1985), with a constant smoothing factor $\alpha \in (0,1)$.

A pseudocode used for this first-stage process is furnished in Appendix C.1.

If the vehicle capacity constraint is non-binding, (16b) unveils that the optimal stop spacing at location $x$ increases with cross-sectional patron flows, $C_{cw}(x)$ and $C_{ccw}(x)$, and decreases with the total OD density, $P_{cw}(x) + Q_{cw}(x) + P_{ccw}(x) + Q_{ccw}(x)$. This is consistent with a previous finding reported for optimal all-stop designs (Wirasinghe and Ghoneim, 1981). The equation also shows that: (i) optimal stop spacing decreases as the number of lines, $m_{cw}$ or $m_{ccw}$, grows; and (ii) spacing between two consecutive non-transfer stops of the same line, which is roughly $\tilde{s}(x)m_{cw}$, increases with the number of lines.[8]

### 3.2 Stage 2: finding optimal scalar values

With $s^*(x_j)$ and $T^*(x_j)$ in hand, $m_{cw}^*$ and $m_{ccw}^*$ are found via enumeration in $\{1,2,3,4\}$, and $H_{cw}^*$ and $H_{ccw}^*$ are found via iteration. The problem is convex with respect to $H_{cw}$ and $H_{ccw}$. This is because (15) is a linear combination of $H_{cw}$, $H_{ccw}$, $H_{cw}^{-1}$ and $H_{ccw}^{-1}$ with non-negative coefficients; and (14c-f) are boundary constraints of $H_{cw}$ and $H_{ccw}$. The $H_{cw}^*$, for example, is thus updated by:

---

[7] We consider $T^*(x_j) \leq 30, \forall j$, to be compatible with what is observed in real settings.

[8] If we omit some minor items in the RHS of (16b) and let $m_{cw} = m_{ccw}$, then (16b) simplifies to: $\tilde{s}(x) \approx \sqrt{\frac{C_{cw}(x)+C_{ccw}(x)}{P_{cw}(x)+Q_{cw}(x)+P_{ccw}(x)+Q_{ccw}(x)}} \cdot \frac{4v_w\tau}{m_{cw}}$. This result reveals that optimal stop spacing (under a symmetric design with $m_{cw} = m_{ccw}$) is roughly proportional to $\sqrt{\frac{C_{cw}(x)+C_{ccw}(x)}{P_{cw}(x)+Q_{cw}(x)+P_{ccw}(x)+Q_{ccw}(x)}}$ and inversely proportional to $\sqrt{m_{cw}}$. It also indicates that single-line stop spacing, $\tilde{s}(x)m_{cw}$, is proportional to $\sqrt{m_{cw}}$.

$$H_{cw}^* = \text{mid}\left\{H_{min} + I(m_{cw} > 1) \cdot \tau, \widetilde{H}_{cw}, \frac{K}{\max\limits_{1 \leq j \leq n}\{C_{cw}(x_j) + \frac{1}{2}T(x_j)s(x_j)(b_{cw}(x_j) + b_{ccw}(x_j))\}}\right\}, \tag{17a}$$

where

$$\widetilde{H}_{cw} = \sqrt{\frac{\frac{\pi_k L}{\mu} + \frac{\pi_h}{\mu}\sum_{j=1}^{n}\left[\frac{1}{v} + \frac{\tau}{s^*(x_j)}\left(\frac{1}{m_{cw}} + \frac{m_{cw}-1}{m_{cw}T^*(x_j)}\right)\right]\Delta x}{\frac{2m_{cw}-1}{2}\Lambda_{cw} + \sum_{j=1}^{n}\left[-\frac{m_{cw}-1}{2}\frac{P_{cw}(x_j) + Q_{cw}(x_j)}{T^*(x_j)} - \frac{m_{cw}}{2}b_{cw}(x_j) + \frac{m_{cw}}{2}b_{ccw}(x_j)\right]\Delta x}}, \tag{17b}$$

and the mid function returns the middle value for the three arguments in (17a). The first of these arguments comes from constraint (14c); the second, $\widetilde{H}_{cw}$, is derived from the first-order condition of $GC$ with respect to $H_{cw}$; and the third is the result of constraint (14e). Similar expressions yield $H_{ccw}^*$.

A pseudocode used for this second-stage process is given in Appendix C.2. The two stages are iterated until they converge.

Result (17b) unveils that the optimal service headway for all lines decreases as the number of lines grows. However, the headway for each line, $m_{cw}\widetilde{H}_{cw}$, obviously increases with $m_{cw}$. This is readily explainable: when more lines are deployed, the overall service frequency (inverse of the headway) increases, while each line's service frequency would decrease because the demand served per line diminishes.

## 4 Model Testing

The two-stage heuristic described in Section 3 cannot guarantee globally-optimized solutions. Optimality gaps between heuristic solutions and lower bounds are therefore evaluated below. Computation times and suitable discretization for the spatial interval are reported as well. A description of the set-up used for this testing comes first.

### 4.1 Experimental set-up

Cases of directionally-symmetric travel demands are alone considered in the interest of brevity. These demands are described by a density function of the form:

$$\lambda(x, y) = p(x)\theta(l(x, y))\Lambda, \tag{18}$$

where: $p(x)$ is the probability density function (PDF) of trip origins over the corridor; $\theta(l(x, y))$ is the PDF of trip length, with $l(x, y) = \min(|x - y|, L - |x - y|)$ denoting trip length, $0 < l \leq L/2$; and $\Lambda = \Lambda_{cw} = \Lambda_{ccw}$ is the total demand in each travel direction. The function (18) enables separate examination of: (i) spatial-heterogeneity in trip origins; and (ii) a distribution of trip lengths.

Demands circulated around a loop corridor, with $L = 40$ km. The $p(x)$ was assumed to be a truncated normal PDF with mean $\frac{L}{2}$ and variance $\sigma_o^2$, truncated by the interval $[0, L]$. A larger value of $\sigma_o$ indicates a "flatter" curve of $p(x)$, which means that trip origins are less varied in space. A $\sigma_o = \infty$ means $p(x)$ is a constant and trip origins are uniformly distributed. The $\theta(l)$ was assumed to be a uniformly-distributed PDF, denoted $\mathcal{U}(E_l - \sqrt{3}\sigma_l, E_l + \sqrt{3}\sigma_l)$; where $E_l$ and $\sigma_l$ are the mean and standard deviation of trip length, respectively. Multiple values of $\sigma_o$, $E_l$ and $\sigma_l$ were used to examine an array of demand patterns.

Input values are shown in Table 1. These include: two values of time used to represent rich and poor cities; and two arrays of average demand density, $\Lambda/L$, one that is suitable for service by bus, the other by rail. In total, 144 cases were examined.

Table 1. Demand parameters and values of time

| Notation | Set of values | Description |
|---|---|---|
| $\sigma_o$ (km) | $\{\infty, 8, 4\}$ | Standard deviation in trip origin distribution. The three values indicate: no-fluctuation, or uniform demand; low-fluctuation; and high-fluctuation, respectively. |
| $E_l$ (km) | $\{8, 12\}$ | Average trip length |
| $\sigma_l$ (km) | $\{2, 4\}$ | Standard deviation in trip length |
| $\mu$ (\$/h) | $\{5, 10\}$ | Value of time |
| $\frac{\Lambda}{L}$ (trips/km/h) | $\{37.5, 75, 150\}$ for bus service; $\{250, 500, 1000\}$ for rail service | Average demand density for each travel direction |

Parameter values for bus and rail services were borrowed from previous studies (Daganzo, 2010; Sivakumaran et al., 2014; Gu et al., 2016), and are shown in Table 2. A low walking speed of $v_w = 2$ km/h was used to account for delays caused by traffic signals. A transfer penalty, $C_t$, was set to 1 min/transfer.

As regards parameters needed by the solution algorithms, the smoothing factor, $\alpha$, was set to 0.5, and tolerances used for convergence (see again Appendix C) to 0.0001. The discretization interval, $\Delta x$, was set to a relatively large value of 0.5 km, following numerical tests of the trade-offs between discretization error and computer runtime.

Table 2. Cost and operational parameters for bus and rail

|  | $\pi_k$ (\$/veh·km) | $\pi_h$ (\$/veh·h) | $\pi_i$ (\$/km/h) | $\pi_s$ (\$/stop/h) |
|---|---|---|---|---|
| Bus | 0.59 | $2.66 + 3\mu$ | $6 + 0.2\mu$ | $0.42 + 0.014\mu$ |
| Rail | 2.20 | $101 + 5\mu$ | $594 + 19.8\mu$ | $294 + 9.8\mu$ |
|  | $\tau$ (sec) | $v$ (km/h) | $K$ (patrons/veh) | $H_{min}$ (min) |
| Bus | 30 | 25 | 80 | 1 |
| Rail | 45 | 60 | 3000 | 1.5 |

## 4.2 Optimality gaps, approximation errors, and computation times

Lower-bound solutions for our two-stage heuristic are presented in Appendix D. Optimality gaps between the heuristic and the resulting lower bounds were assessed for the 144 cases cited in Section 4.1. The average gap across these cases was 0.7%; and the maximum of these gaps was 2.8%. These small differences reflect well on the heuristic, and indicate that backtracking costs are small. The latter were ignored in the lower bounds, as is made clear in Appendix D.[9]

To examine the accuracy of our CA model, we convert the heuristic solution into exact locations of stops and skip-stop bays using an algorithm described in Appendix E,[10] and re-calculate the user and agency costs using an algorithm presented in Appendix F.1. Percentage errors were then computed between the recalculated cost components and those of the heuristic solution. Outcomes were again favorable, as described in Appendix F.2.

As an important aside, average computation time for the 144 cases was 151s. Approximately 90% of the times were devoted to preprocessing demands, as per (1a)-(3). Less than 20s was typically needed to find solutions via the heuristic, indicating that the algorithm is efficient. Generating a stop location plan generally required only a second or two.

---

[9] Backtracking was not ignored in the heuristic, however, since this would have led to covering an entire corridor with only a single skip-stop bay.

[10] The algorithm ensures that the exact stop spacings and (integer-valued) numbers of stops in each skip-stop bay respectively match the $s^*(x_j)$ and $T^*(x_j)$ obtained for each $x_j$ by the two-stage heuristic (see again Section 3.1).

### 4.3 Deployment issues

Tests further show that $s^*(x)$ and $T^*(x)$ closely match exact locations of transfer and non-transfer stops, as generated using the algorithm in Appendix E. Examples are shown in Figures 4a and b.[11] Interested readers may also refer to Mei (2019) for a real-world case study.

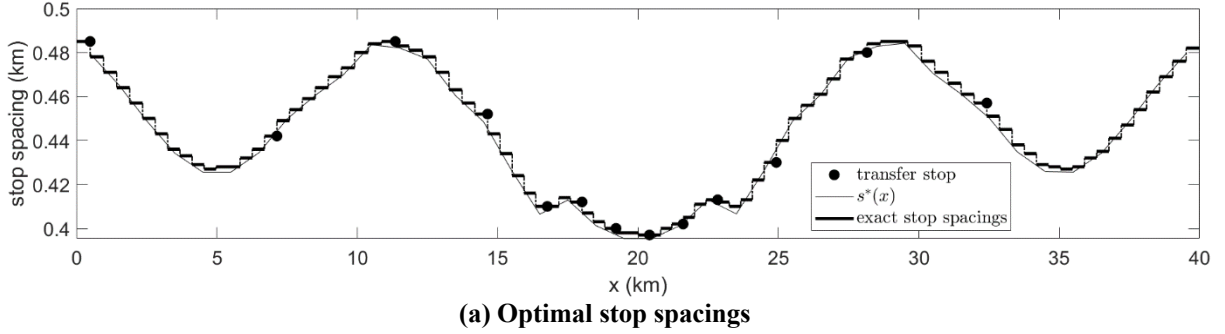

(a) Optimal stop spacings

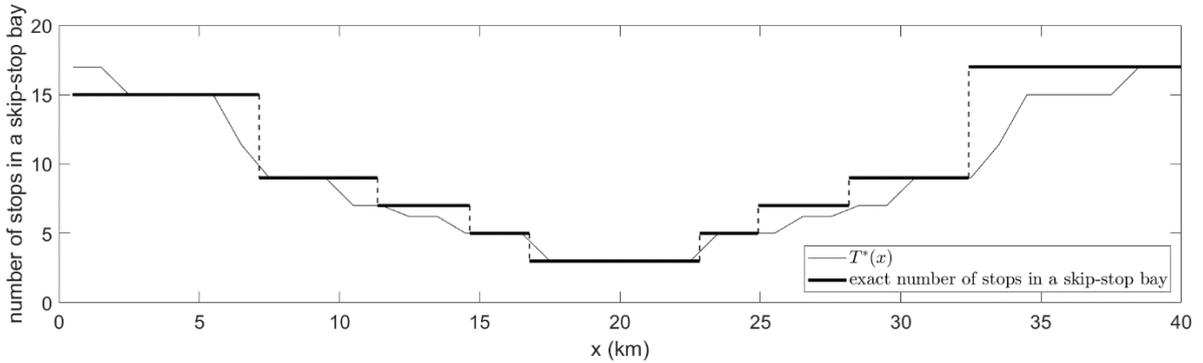

(b) Optimal numbers of stops in each skip-stop bay
Figure 4. Optimal AB-type design for an example bus corridor

## 5 Numerical Analysis

Costs saved via optimized AB-type service relative to optimized all-stop service are compared below. Parametric analysis of corridors served by buses are examined in Section 5.1, and rail corridors in Section 5.2. Key inputs (transfer penalties and patron access speeds) are studied parametrically in Section 5.3.

### 5.1 Bus corridors

Consider for illustration high-wage cities with $\mu = \$20/h$ that are characterized by moderately heterogeneous demands with $\sigma_o = 8$ km. The $GC$s saved by AB-type service is shown in Figure 5a as functions of demand density, and for various trip-length distributions. Note how skip-stop service reduces $GC$, except for cases characterized by both low-demand densities and short-length trips. Further note how cost savings grow with growing demand densities and trip lengths. Since skip-stop service tends to increase agency costs and reduce user costs, its advantage grows as the system serves greater patron-kms of travel.[12]

Effects of $\mu$ and $\sigma_o$ are made evident in Figure 5b. In these cases, the $E_l$ and $\sigma_l$ are fixed at 8km and 4km, respectively. The missing data reflect conditions in which the buses' patron-carrying capacity

---

[11] The curves were derived for: $E_l = 8$ km; $\sigma_l = 2$ km; $\sigma_o = 4$ km; $\frac{\Lambda}{L} = 37.5$ trips/km/h; $\mu = 20$ \$/h; and 2 lines in each travel direction.

[12] Analysis not shown here (for brevity) indicates that skip-stop service can even save agency cost when demands are highly heterogeneous over space. In these instances, skip-stop service requires shorter round-trip times, and thus fewer buses.

($K = 80$ patrons/bus, as specified in Table 2) was exceeded. Still, consideration of the curves reveals that, all else equal, cost saved by skip-stop service grows with: greater spatial-heterogeneity in demand; and higher values of time.

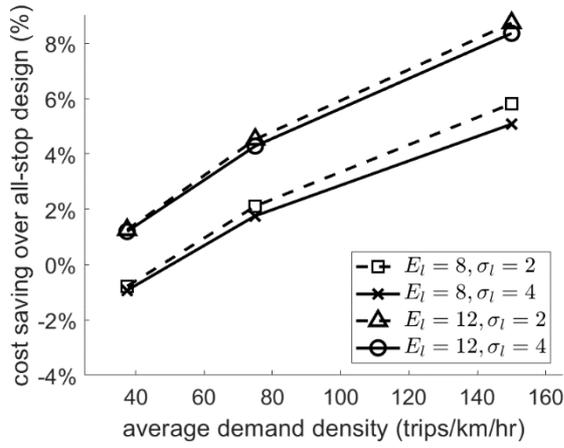 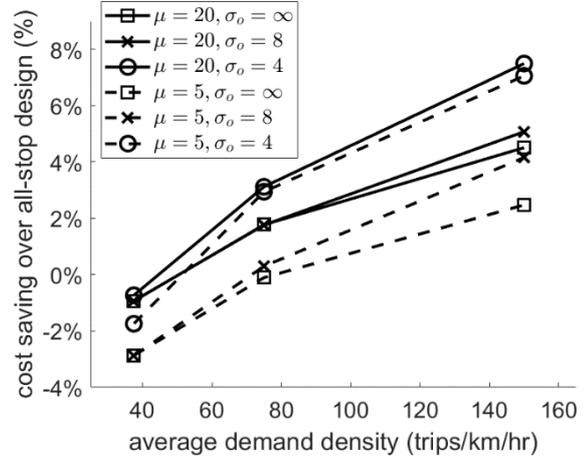

(a) Demand with less-fluctuating trip origins in high-wage cities ($\sigma_o = 8$ km, $\mu = 20$ $/h)

(b) Demand with short and highly-varied trip lengths ($E_l = 8$ km, $\sigma_l = 4$ km)

Figure 5. Percentage generalized cost savings of AB-type bus systems

### 5.2 Rail corridors

We find that skip-stop service almost always saves $GC$ in high-demand corridors served by rail. This is evident in Figures 6a and b. The trends in these curves are similar to those for bus corridors in the previous section.

### 5.3 Transfer penalty, access speed, and weight for backtracking cost

Suitable choices for transfer penalty, $C_t$, and patron access speed, $v_w$, are liable to be case-specific, and are therefore examined parametrically.[13] We do so here for rail corridors with $\sigma_o = 8$ km, $E_l = 12$ km, $\sigma_l = 4$ km, and $\mu = \$20$/h. Note from Figure 7a how costs saved by skip-stop service diminish under higher transfer penalties. Savings also diminish as access speed increases, as evident in Figure 7b. Savings nearly vanish under high access speed of 8 km/h. The higher speed justifies larger stop spacings, with fewer stops therefore visited each trip. This makes all-stop service a better option.

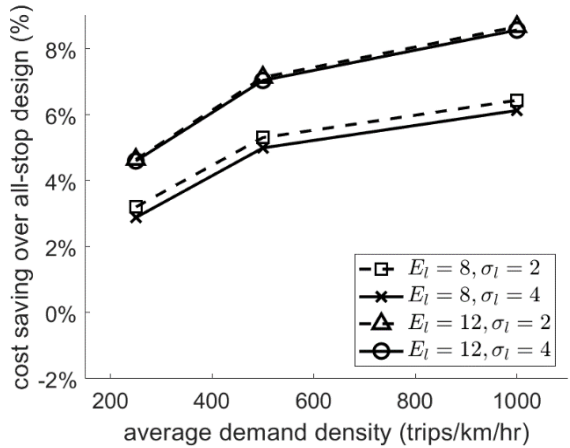 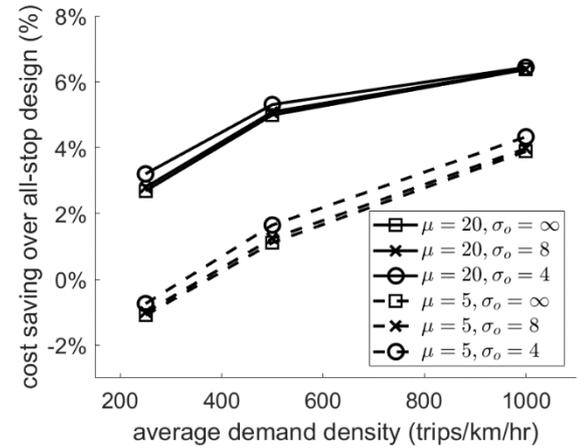

(a) Demand with highly-fluctuating trip origins in high-wage cities ($\sigma_o = 4$ km, $\mu = 20$ $/h)

(b) Demand with short and highly-varied trip lengths ($E_l = 8$ km, $\sigma_l = 4$ km)

---

[13] Values of $C_t$ were varied from 1-2 minutes as per Guo and Wilson (2004) and Fan et al. (2018); and $v_w$ varied from 2-8 km/h, as per Wu et al. (2020).

**Figure 6. Percentage generalized cost savings of AB-type rail systems**

In addition, to account for patrons' aversion to backtracking, we multiply the backtracking cost in $GC$ by a weight, which is assumed to vary from 1 to 3. Cost savings are plotted against this weight in Figure 7c for a bus corridor with $\sigma_o = 4$ km, $E_l = 8$ km, $\sigma_l = 2$ km, $\frac{\Lambda}{L} = 75$ trips/km/h, and $\mu = \$20$/h. The figure shows that the benefit of skip-stop design is insensitive to the cost specified for backtracking.

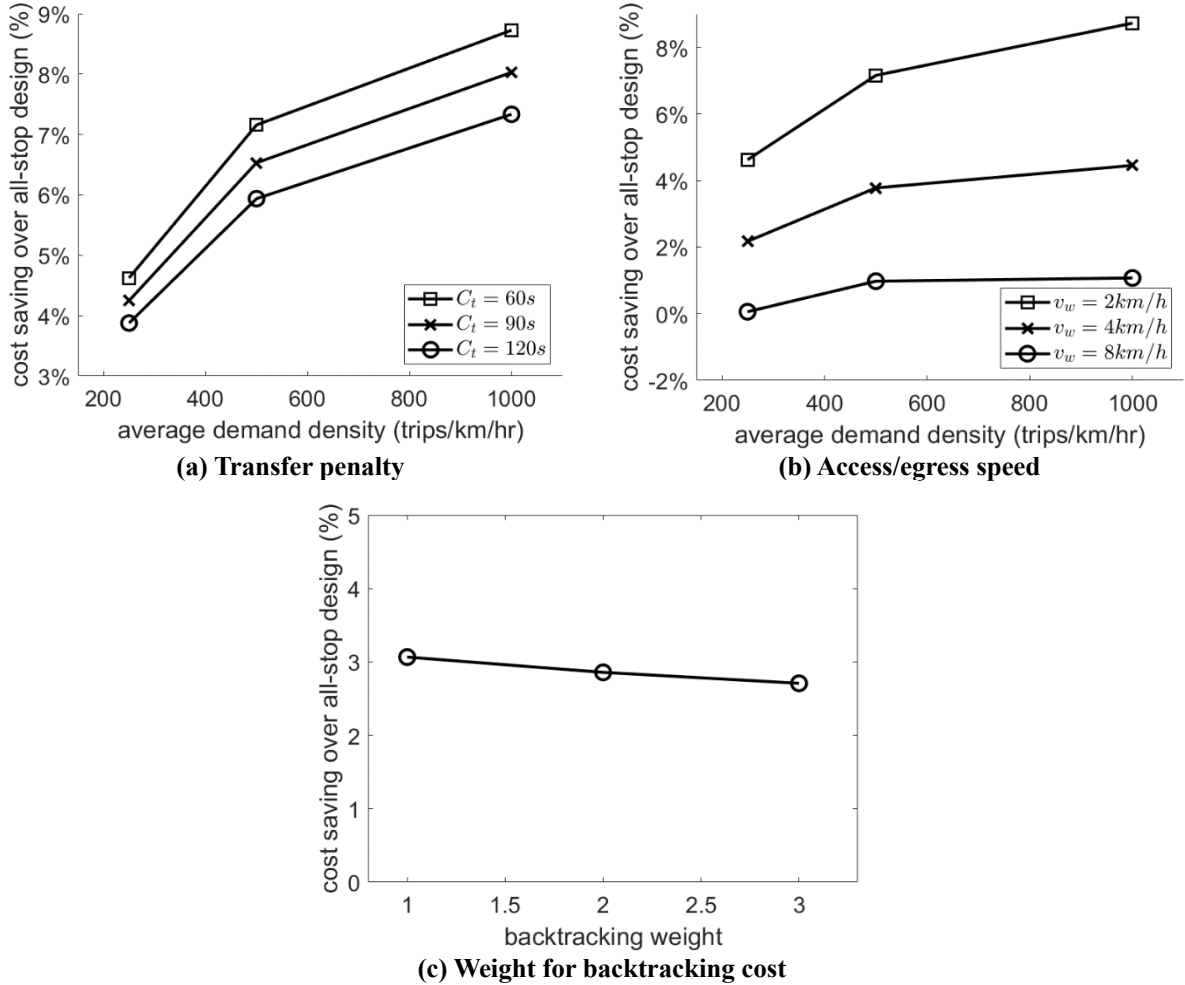

(a) Transfer penalty  (b) Access/egress speed

(c) Weight for backtracking cost

**Figure 7. Sensitivity of cost savings to transfer penalty, access/egress speed, and backtracking cost**

## 6 Conclusions

The present work has explored a particular form of skip-stop service, termed AB-type service, in which non-transfer stops are arranged in certain order along distinct lines. The continuous-approximation model developed for this exploration can accommodate demand patterns that vary along the length of the corridor. The model furnishes stop densities and route plans (a corridor's number of lines and the stops visited each run) that may vary with location, and by direction. Calculus of variations and direct search methods were used to formulate a heuristic solution method. Costs compared against a lower bound indicate that the heuristic produces near-optimal designs with short runtimes. These solutions can identify exact stop locations using the algorithm in Appendix E. The resulting designs can be further fine-tuned to meet physical constraints in real settings; e.g., to avoid placing stops at intersections.

Further comparisons against optimized all-stop designs confirm that AB-type service is often the lower-cost option, particularly on high-demand corridors traveled by wealthier patrons with high values of time. In these cases, AB-type service produced cost savings up to 8%.

These findings are consistent with earlier ones in Gu et al. (2016). The model developed in that work, however, assumed that travel demands are spatially homogeneous. As a result, the earlier model produces stop densities and route plans that are spatially uniform. Thanks to its accommodation of spatially-varying demands, the present model unveils new findings in addition to those noted above. We now find, for example, that optimal stop spacing and skip-stop bay length are both inversely correlated to the local OD density, and that stop spacing is positively correlated to the local cross-sectional patron flows. In addition, we confirm that AB-type service is especially beneficial when trip origins are spatially heterogeneous.

## Acknowledgements

The research was supported by General Research Funds (No. 15217415 and No. 15224818) provided by the Research Grants Council of Hong Kong, and funds provided by the National Natural Science Foundation of China (No. 51608455) and Sichuan Provincial Science & Technology Innovation Cooperation Funds (No. 2020YFH0038).

## Appendix A. Table of notations

Table A1. List of variables, functions, and parameters

| Notation | Description | Notation | Description |
|---|---|---|---|
| \multicolumn{4}{c}{*Decision variables and functions*} ||||
| $m_{cw}$ | Number of CW skip-stop lines | $H_{ccw}$ | Headway in the CCW direction (h) |
| $m_{ccw}$ | Number of CCW skip-stop lines | $s(x)$ | Stop spacing at $x$ (km) |
| $H_{cw}$ | Headway in the CW direction (h) | $T(x)$ | Number of stops in a skip-stop bay at $x$ |
| \multicolumn{4}{c}{*Demand variables and functions*} ||||
| $\lambda(x,y)$ | Demand density from origin $x$ to destination $y$ (trip/km²/h) | $Q_{cw}(x)$ | CW trip destination density at $x$ (trip/km/h) |
| | | $Q_{ccw}(x)$ | CCW trip destination density at $x$ (trip/km/h) |
| $\Lambda_{cw}$ | Total CW demand (trip/h) | $l$ | Trip length (km) |
| $\Lambda_{ccw}$ | Total CCW demand (trip/h) | $p(x)$ | PDF of trip origins (under symmetric demand) |
| $C_{cw}(x)$ | On-board flow of CW trips at $x$ (trip/h) | $\theta(l)$ | PDF of trip lengths (under symmetric demand) |
| $C_{ccw}(x)$ | On-board flow of CCW trips at $x$ (trip/h) | $\sigma_o$ | Standard deviation of origin distribution (km) |
| $P_{cw}(x)$ | CW trip origin density at $x$ (trip/km/h) | $E_l$ | Mean trip length (km) |
| $P_{ccw}(x)$ | CCW trip origin density at $x$ (trip/km/h) | $\sigma_l$ | Standard deviation of trip length (km) |
| \multicolumn{4}{c}{*Cost terms and parameters*} ||||
| $\pi_k$ | Unit distance-based operating cost ($/veh·km) | $\pi_i$ | Amortized unit cost of line infrastructure ($/km/h) |
| $\pi_h$ | Unit time-based operating cost ($/veh·h) | $\pi_s$ | Amortized unit cost per stop ($/stop/h) |
| $C_t$ | Unit transfer penalty cost (h) | $AC_K$ | Distance-based operating cost (h/h) |
| $UT_a$ | Patrons' total access and egress time (h/h) | $AC_H$ | Time-based operating cost (h/h) |
| $UT_w$ | Patrons' total wait time (h/h) | $AC_I$ | Line infrastructure cost (h/h) |
| $UT_v$ | Patrons' total in-vehicle travel time (h/h) | $AC_S$ | Stop infrastructure cost (h/h) |
| $UT_t$ | Total transfer penalty (h/h) | $GC$ | Generalized cost (h/h) |
| \multicolumn{4}{c}{*Other parameters and variables*} ||||
| $b_{cw}(x)$ | Density of backtracking trips at $x$ in the CW direction (trip/km/h) | $\bar{\lambda}_{cw}(x)$ | Average density of CW trips contained in a skip-stop bay at $x$ (trip/km/h) |
| $b_{ccw}(x)$ | Density of backtracking trips at $x$ in the CCW direction (trip/km/h) | $\bar{\lambda}_{ccw}(x)$ | Average density of CCW trips contained in a skip-stop bay at $x$ (trip/km/h) |
| $v_{cw}(x)$ | Commercial speed at $x$ in the CW direction | $\mu$ | Patrons' value of time ($/h) |
| $v_{ccw}(x)$ | Commercial speed at $x$ in the CCW direction | $H_{min}$ | Minimum headway (h) |
| $L$ | Corridor length (km) | $K$ | Vehicle capacity (patron/veh) |
| $v_w$ | Walking speed (km/h) | $\tau$ | Vehicle dwell time per stop (h) |
| $v$ | Vehicle cruise speed (km/h) | | |

# Appendix B. Formulas related to user cost metrics

## B.1 Average density of backtracking trips in (5)

Denote $\bar{\lambda}_{cw}(x)$ (trips/km/h) as the average density of CW "contained" trips at $x$, where a contained trip is one in which origin and destination are located in the same skip-stop bay. The $\bar{\lambda}_{cw}(x)$ is calculated by dividing the total number of contained trips in the skip-stop bay containing $x$ by the length of that bay, i.e.,

$$\bar{\lambda}_{cw}(x) = \frac{\int_{z=U(x)}^{D(x)} \int_{y=z}^{D(x)} \lambda(z,y) dz dy}{D(x)-U(x)} \approx \int_{z=x-\frac{T(x)s(x)}{2}}^{x+\frac{T(x)s(x)}{2}} \int_{y=z}^{x+\frac{T(x)s(x)}{2}} \lambda(z,y) dz dy \ / T(x)s(x) \qquad (B1)$$

where $U(x)$ and $D(x)$ denote the locations of the upstream and downstream transfer stops that bound the bay. Here they are approximated by $x - T(x)s(x)/2$ and $x + T(x)s(x)/2$, respectively.

The probability that a contained trip at $x$ is a backtracking trip is $\frac{T(x)-1}{T(x)} \cdot \frac{(m_{cw}-1)(T(x)-1)}{m_{cw}T(x)}$. The first fraction in this formula is the probability that a trip's origin stop is a non-transfer stop in the bay; and the second fraction is the probability that the trip's destination stop is a non-transfer one in the same bay on a different line. We then have $b_{cw}(x) = \frac{(m_{cw}-1)(T(x)-1)^2}{m_{cw}T(x)^2} \bar{\lambda}_{cw}(x)$, which yields (5).

## B.2 Total wait time in (6)

The average wait times and hourly trip numbers for each of the five trip types in CW direction are summarized in the following table.

**Table B1. Wait times for five CW trip types**

| Trip type | Average wait time | Number of trips per hour |
|---|---|---|
| 1 | $W_1 = H_{cw}/2$ | $N_1 = \iint_{D_{cw}} \frac{\lambda(x,y)}{T(x)T(y)} dxdy$ |
| 2 | $W_2 = m_{cw}H_{cw}/2$ | $N_2 = \iint_{D_{cw}} \frac{(T(x)+T(y)-2)\lambda(x,y)}{T(x)T(y)} dxdy$ |
| 3 | $W_3 = m_{cw}H_{cw}/2$ | $N_3 = \iint_{D_{cw}} \frac{(T(x)-1)(T(y)-1)\lambda(x,y)}{m_{cw}T(x)T(y)} dxdy$ |
| 4 | $W_4 = (m_{cw}H_{cw} + m_{ccw}H_{ccw})/2$ | $N_4 = \int_{x=0}^{L} b_{cw}(x) dx$ |
| 5 | $W_5 = m_{cw}H_{cw}$ | $N_5 = \Lambda_{cw} - \sum_{i=1}^{4} N_i$ |

The total wait time for CW trips is then obtained by plugging the formulas of $N_i$ and $W_i$ ($i = 1,2,3,4,5$) into $\sum_{i=1}^{5} N_i W_i$:

$$\sum_{i=1}^{5} N_i W_i = \iint_{D_{cw}} \left( \frac{1}{T(x)T(y)} \frac{H_{cw}}{2} + \frac{(T(x)+T(y)-2)}{T(x)T(y)} \frac{m_{cw}H_{cw}}{2} + \frac{(T(x)-1)(T(y)-1)}{m_{cw}T(x)T(y)} \frac{m_{cw}H_{cw}}{2} \right) \lambda(x,y) dxdy +$$
$$N_4 \frac{m_{cw}H_{cw}+m_{ccw}H_{ccw}}{2} + \left( \iint_{D_{cw}} \left( \frac{m_{cw}-1}{m_{cw}} \frac{(T(x)-1)(T(y)-1)}{T(x)T(y)} \right) m_{cw}H_{cw} \lambda(x,y) dxdy - N_4 \cdot m_{cw}H_{cw} \right)$$
$$= \frac{H_{cw}}{2} \iint_{D_{cw}} \lambda(x,y) \frac{1+m_{cw}T(x)+m_{cw}T(y)-2m_{cw}+(2m_{cw}-1)(T(x)T(y)-T(x)-T(y)+1)}{T(x)T(y)} dxdy +$$
$$\frac{m_{ccw}H_{ccw}-m_{cw}H_{cw}}{2} N_4$$
$$= \frac{H_{cw}}{2} \iint_{D_{cw}} \lambda(x,y) \left( 2m_{cw} - 1 - \frac{m_{cw}-1}{T(x)} - \frac{m_{cw}-1}{T(y)} \right) dxdy + \frac{m_{ccw}H_{ccw}-m_{cw}H_{cw}}{2} N_4$$
$$= \frac{(2m_{cw}-1)H_{cw}}{2} \iint_{D_{cw}} \lambda(x,y) dxdy - \frac{(m_{cw}-1)H_{cw}}{2} \left( \iint_{D_{cw}} \frac{\lambda(x,y)}{T(x)} dydx + \iint_{D_{cw}} \frac{\lambda(x,y)}{T(y)} dxdy \right) +$$
$$\frac{m_{ccw}H_{ccw}-m_{cw}H_{cw}}{2} N_4$$
$$= \frac{(2m_{cw}-1)H_{cw}\Lambda_{cw}}{2} - \frac{(m_{cw}-1)H_{cw}}{2} \int_{x=0}^{L} \int_{y=x}^{x+\frac{L}{2}} \frac{\lambda(x,y)}{T(x)} dydx - \int_{y=0}^{L} \int_{x=y-\frac{L}{2}}^{y} \frac{\lambda(x,y)}{T(y)} dxdy +$$
$$\frac{m_{ccw}H_{ccw}-m_{cw}H_{cw}}{2} N_4$$

$$= \frac{(2m_{cw}-1)H_{cw}\Lambda_{cw}}{2} - \frac{(m_{cw}-1)H_{cw}}{2}\left(\int_{x=0}^{L}\frac{P_{cw}(x)}{T(x)}dx + \int_{y=0}^{L}\frac{Q_{cw}(y)}{T(y)}dy\right) + \frac{m_{ccw}H_{ccw}-m_{cw}H_{cw}}{2}N_4$$

$$= \frac{(2m_{cw}-1)H_{cw}\Lambda_{cw}}{2} - \frac{(m_{cw}-1)H_{cw}}{2}\int_{x=0}^{L}\frac{P_{cw}(x)+Q_{cw}(x)}{T(x)}dx + \frac{m_{ccw}H_{ccw}-m_{cw}H_{cw}}{2}\int_{x=0}^{L}b_{cw}(x)dx \quad (B2)$$

The fifth equality in (B2) uses the definition of $\Lambda_{cw}$, and the sixth equality uses the definitions of $P_{cw}(x)$ and $Q_{cw}(y)$.

The patron wait time for all CCW trips can be obtained by swapping subscripts "$cw$" with "$ccw$" in (B2). Combining both directions, we have (6).

### B.3 Approximation of the transfer penalty in (9)

$$\iint_{D_{cw}} \lambda(x,y)\frac{m_{cw}-1}{m_{cw}}\left(1-\frac{1}{T(x)}\right)\left(1-\frac{1}{T(y)}\right)dxdy$$

$$= \frac{m_{cw}-1}{m_{cw}}\left(\iint_{D_{cw}}\lambda(x,y)dxdy - \iint_{D_{cw}}\frac{\lambda(x,y)}{T(x)}dydx - \iint_{D_{cw}}\frac{\lambda(x,y)}{T(y)}dxdy + \iint_{D_{cw}}\frac{\lambda(x,y)}{T(x)T(y)}dxdy\right)$$

$$= \frac{m_{cw}-1}{m_{cw}}\left(\Lambda_{cw} - \int_{x=0}^{L}\frac{P_{cw}(x)+Q_{cw}(x)}{T(x)}dx + \iint_{D_{cw}}\frac{\lambda(x,y)}{T(x)T(y)}dxdy\right)$$

$$\approx \frac{m_{cw}-1}{m_{cw}}\left(\Lambda_{cw} - \int_{x=0}^{L}\frac{P_{cw}(x)+Q_{cw}(x)}{T(x)}dx + \frac{1}{2}\iint_{D_{cw}}\left(\frac{\lambda(x,y)}{T^2(x)}+\frac{\lambda(x,y)}{T^2(y)}\right)dxdy\right)$$

$$= \frac{m_{cw}-1}{m_{cw}}\left(\Lambda_{cw} - \int_{x=0}^{L}\frac{P_{cw}(x)+Q_{cw}(x)}{T(x)}dx + \frac{1}{2}\int_{x=0}^{L}\frac{P_{cw}(x)+Q_{cw}(x)}{T^2(x)}dx\right)$$

$$= \frac{m_{cw}-1}{m_{cw}}\left(\Lambda_{cw} - \int_{x=0}^{L}\frac{(P_{cw}(x)+Q_{cw}(x))(2T(x)-1)}{2T^2(x)}dx\right) \quad (B3)$$

The third and fifth lines in (B3) were derived using the definitions of $\Lambda_{cw}$, $P_{cw}(x)$ and $Q_{cw}(x)$. The approximation in the fourth line was used to reduce the last double integral term to a single one. Since $\frac{\lambda(x,y)}{T(x)T(y)} \leq \frac{1}{2}\left(\frac{\lambda(x,y)}{T^2(x)}+\frac{\lambda(x,y)}{T^2(y)}\right)$, the above approximation is conservative. The number of transfers for CCW trips is obtained by replacing subscript "$cw$" in the above result by "$ccw$". Equation (9) then follows.

## Appendix C. Pseudocodes of the solution algorithms

### C.1 Pseudocode of the first-stage solution process

**Algorithm 1: Finding $s^*(x_j)$ and $T^*(x_j)$ for $j = 1, \ldots, n$, given $m_{cw}, m_{ccw}, H_{cw}$ and $H_{ccw}$.**

For $j = 1, 2, \ldots, n$:
    Let $x_j = (j - 0.5)\Delta x$;
    Initialize $b_{cw}(x_j) = b_{ccw}(x_j) = 0$;
    Do:
        For each $T(x_j)$ in the candidate set $\{1,2,\ldots,30\}$:
            Update $s^*(x_j)$ using Equations (16a-b);
            Record the lowest $G(\cdot, x_j)$ by far, and the associated $T^*(x_j)$ and $s^*(x_j)$;
        End For
        Calculate $\hat{b}_{cw}(x_j)$ using (5), and calculate $\hat{b}_{ccw}(x_j)$ similarly;
        $b_{cw}(x_j) \leftarrow (1-\alpha)b_{cw}(x_j) + \alpha\hat{b}_{cw}(x_j)$;
        $b_{ccw}(x_j) \leftarrow (1-\alpha)b_{ccw}(x_j) + \alpha\hat{b}_{ccw}(x_j)$;
    Until $|b_{cw}(x_j) - \hat{b}_{cw}(x_j)| + |b_{ccw}(x_j) - \hat{b}_{ccw}(x_j)| \leq \epsilon_1$
End For
**Output** $s^*(x_j)$ and $T^*(x_j)$ for $j = 1, 2, \ldots, n$.

Here $\epsilon_1$ is a convergence tolerance.

### C.2 Pseudocode of the second-stage solution process

**Algorithm 2: Finding $m_{cw}^*, m_{ccw}^*, H_{cw}^*$ and $H_{ccw}^*$.**

For each pair of $(m_{cw}, m_{ccw})$ in $\{1,2,3,4\} \times \{1,2,3,4\}$:
    Initialize $H_{cw}$ and $H_{ccw}$ that satisfy (14c-f) by letting $b_{cw}(x) = b_{ccw}(x) = 0$;
    Do:
        Find $s^*(x_j)$ and $T^*(x_j)$ for all $j = 1,2,\ldots,n$ using **Algorithm 1**;
        Calculate $H_{cw}^*$ using Equation (17a-b), and calculate $H_{ccw}^*$ similarly;
    Until $|H_{cw} - H_{cw}^*| + |H_{ccw} - H_{ccw}^*| \leq \epsilon_2$
    Record the $(m_{cw}^*, m_{ccw}^*)$ that yields the lowest $GC$, and the associated $H_{cw}^*, H_{ccw}^*, s^*(x_j)$ and $T^*(x_j)$;
End For
**Output** $m_{cw}^*, m_{ccw}^*, H_{cw}^*, H_{ccw}^*, s^*(x_j | j = 1,2,\ldots,n), T^*(x_j | j = 1,2,\ldots,n)$.

Here $\epsilon_2$ is another convergence tolerance.

**Appendix D. Lower bound of (14a-g)**

Constraints (14e-f) are modified by ignoring $T(x)s(x)(b_{cw}(x) + b_{ccw}(x))/2$ in the denominator of their RHS. This results in looser constraints, and thus a greater feasible solution region for the relaxed program.

Reduce the original objective function (14a) in two steps as explained below.

Step 1. Remove the cost terms related to $b_{cw}(x)$ and $b_{ccw}(x)$, including parts of the RHS of (6) and (8). These removed terms sum to:

$$\int_{x=0}^{L} \left[ \frac{m_{ccw}H_{ccw} - m_{cw}H_{cw}}{2}(b_{cw}(x) - b_{ccw}(x)) + \left[\frac{s(x)T(x)}{3v} + \frac{\tau}{6}\left(\frac{T(x)-1}{m_{cw}} + \frac{T(x)-1}{m_{ccw}} + 2\right)\right](b_{cw}(x) + b_{ccw}(x)) \right] dx. \quad (D1)$$

To guarantee that a lower bound is achieved, (D1) must be non-negative. This is evident under symmetric demand (i.e., $\lambda(x,y) = \lambda(y,x)$ for any $(x,y)$), because in that case, $m_{cw} = m_{ccw}$ and $H_{cw} = H_{ccw}$, and the first term in the integrand is zero while the second term is positive. For asymmetric demand, however, that first term can be negative, and thus theoretically (D1) can be negative too. Therefore, to be exact, the "lower bound" derived in this section is a lower bound only under symmetric demand. Nevertheless, the first term in the integrand is generally close to zero because both $|m_{ccw}H_{ccw} - m_{cw}H_{cw}|$ and $|b_{cw}(x) - b_{ccw}(x)|$ are quite small. Thus, (D1) is still likely to be non-negative under asymmetric demands.

Step 2. Replace $-\frac{2T(x)-1}{2T^2(x)}$ in transfer penalty cost (9) by $-\frac{1}{T(x)}$, since $-\frac{2T(x)-1}{2T^2(x)} > -\frac{1}{T(x)}$.

The new objective function can then be rearranged as:

$$GC_{LB} = \theta + \int_{x=0}^{L} \left( f(x) + \frac{\beta(x)}{T(x)} \right) dx, \quad (D2a)$$

where $\theta$ is a function of scalar variables $H_{cw}, H_{ccw}, m_{cw}$ and $m_{ccw}$ only; and $f(x)$ and $\beta(x)$ are related to $s(x)$ but not to $T(x)$. These parameters are as follows:

$$\theta = \frac{(2m_{cw}-1)\Lambda_{cw}}{2}H_{cw} + \frac{(2m_{ccw}-1)\Lambda_{ccw}}{2}H_{ccw} + C_t \frac{m_{cw}-1}{m_{cw}}\Lambda_{cw} + C_t \frac{m_{ccw}-1}{m_{ccw}}\Lambda_{ccw} + \frac{\pi_k L}{\mu}\left(\frac{1}{H_{cw}} + \frac{1}{H_{ccw}}\right) + \frac{2\pi_i L}{\mu}$$

$$f(x) = \frac{s(x)}{4v_w}\left(P_{cw}(x) + Q_{cw}(x) + P_{ccw}(x) + Q_{ccw}(x)\right) + \left(C_{cw}(x) + \frac{\pi_h}{\mu H_{cw}}\right)\left(\frac{1}{v} + \frac{\tau}{m_{cw}s(x)}\right) + \left(C_{ccw}(x) + \frac{\pi_h}{\mu H_{ccw}}\right)\left(\frac{1}{v} + \frac{\tau}{m_{ccw}s(x)}\right) + \frac{\pi_s}{\mu}\frac{1}{s(x)} \quad (D2b)$$

$$\beta(x) = \frac{\left(C_{cw}(x)+\frac{\pi_h}{\mu H_{cw}}\right)\tau}{s(x)} \frac{m_{cw}-1}{m_{cw}} - \left(\frac{m_{cw}-1}{m_{cw}} C_t + \frac{(m_{cw}-1)H_{cw}}{2}\right)\left(P_{cw}(x)+Q_{cw}(x)\right) +$$
$$\frac{\left(C_{ccw}(x)+\frac{\pi_h}{\mu H_{ccw}}\right)\tau}{s(x)} \frac{m_{ccw}-1}{m_{ccw}} - \left(\frac{m_{ccw}-1}{m_{ccw}} C_t + \frac{(m_{ccw}-1)H_{ccw}}{2}\right)\left(P_{ccw}(x)+Q_{ccw}(x)\right). \tag{D2c}$$

The lower bound problem is then formulated as:

$$\min_{m_{cw},m_{ccw},H_{cw},H_{ccw},s(x),T(x)} GC_{LB} = \theta + \int_{x=0}^{L}\left(f(x)+\frac{\beta(x)}{T(x)}\right)dx \tag{D3a}$$

subject to:
$$m_{cw}, m_{ccw} \in \{1,2,3,4\} \tag{D3b}$$
$$H_{min} + I(m_{cw} > 1)\cdot \tau \leq H_{cw} \leq \frac{K}{\max_{0<x\leq L}\{C_{cw}(x)\}} \tag{D3c}$$
$$H_{min} + I(m_{ccw} > 1)\cdot \tau \leq H_{ccw} \leq \frac{K}{\max_{0<x\leq L}\{C_{ccw}(x)\}} \tag{D3d}$$
$$s(x) > 0, T(x) \geq 1 \tag{D3e}$$

The optimal solution to (D3a-e) is a lower bound to the solution of (14a-g) because the former has a lower objective function value and relaxed constraints. We next solve (D3a-e) via a bi-level method.

At the lower level, fix $m_{cw}, m_{ccw}, H_{cw}$ and $H_{ccw}$, and minimize the integrand $f(x) + \frac{\beta(x)}{T(x)}$ for each $x \in [0,L)$. Note that for any fixed $s(x)$, $f(x) + \frac{\beta(x)}{T(x)}$ is minimized at either $T(x) = 1$ (if $\beta(x)$ is negative) or $T(x) = \infty$ (if $\beta(x)$ is non-negative). Therefore,

$$\min_{s(x),T(x)} f(x) + \frac{\beta(x)}{T(x)} = \min\left\{\min_{s(x)} f(x), \min_{s(x)} f(x) + \beta(x)\right\} \tag{D4}$$

Note further that $f(x)$ and $\beta(x)$ are both convex with respect to $s(x)$. Hence the global optimal solution to the lower-level problem can be obtained using a gradient search method.

At the upper level, optimize $m_{cw}, m_{ccw}, H_{cw}$ and $H_{ccw}$ via exhaustive search. Enumerate $m_{cw}$ and $m_{ccw}$ from $\{1,2,3,4\} \times \{1,2,3,4\}$, and $H_{cw}$ and $H_{ccw}$ from the ranges specified by (D3c) and (D3d) with an interval of 0.1 min. The lowest-cost solution produced by the above procedure is a lower bound to the optimal cost from (14a-g).

## Appendix E. Algorithm for generating stop locations

A 3-step algorithm is proposed below to generate stop locations and route plans using the optimal $s^*(x_j)$ and $T^*(x_j)$ for each $x_j, j = 1,2,...,n$; see Section 3.1. Step 1 generates the stop locations (including both non-transfer and transfer stops) that closely match $s^*(x_j)$ using a method similar to Wirasinghe and Ghoneim (1981). Step 2 selects the transfer stops from the set of stops generated in Step 1 to ensure that the number of stops in each skip-stop bay matches the average $T(x)$ in that bay as close as possible. Step 3 assigns the non-transfer stops to each line. These steps are detailed as follows.

Step 1. Apply spline curve fitting method to fit the continuous function $s(x)$ using $s^*(x_j)$ ($j = 1,2,...,n$). Next, place one stop at every $x$ where $\int_{z=0}^{x} \frac{dz}{s(z)}$ is an integer (the first stop is located at $x = 0$). The resulting stop location set is denoted as $\Omega = \{x_i^S : i = 1,2,...,N^S\}$, where $N^S$ is the number of stops, and $x_i^S$ is the location of the $i$-th stop satisfying $0 = x_1^S < x_2^S < \cdots < x_{N^S}^S < L$. Finally, since $\int_{z=0}^{L} \frac{dz}{s(z)}$ may not be an integer, remove the last stop from $\Omega$ if it is too close to the first one (i.e. if $L - x_{N^S}^S < s(L)/2$).

Step 2. Denote the number of transfer stops as $N^T$, and the stop index of the $k$-th transfer stop ($k = 1, \ldots, N^T$) in $\Omega$ as $u_k$ satisfying $u_1 < u_2 < \cdots < u_{N^T}$. Further denote the $k$-th transfer stop's location as $x_k^T \equiv x_{u_k}^S \in \Omega$. The $N^T$ and $u_k$ ($k = 1, \ldots, N^T$) are generated as follows. Set $u_1 = 1$ (and thus $x_1^T \equiv x_1^S = 0$). Select the other transfer stops recursively; i.e., if $u_k$ and $x_k^T$ are known, find $u_{k+1}$ such that:

$$\min_{u_{k+1} \in \{i=1,2,\ldots,N^S\}} \left| (u_{k+1} - u_k) - \int_{z=x_k^T}^{x_{k+1}^T} T(z) dz / (x_{k+1}^T - x_k^T) \right| \tag{E1a}$$

subject to:
$$x_{k+1}^T = x_{u_{k+1}}^S; \tag{E1b}$$

$$u_{k+1} > u_k; \tag{E1c}$$

$u_{k+1} - u_k - 1$ is an integer multiple of both $m_{cw}$ and $m_{ccw}$. (E1d)

Here $u_{k+1} - u_k$ is the number of stops in the $k$-th skip-stop bay, hence (E1d) is required; and $\int_{z=x_k^T}^{x_{k+1}^T} T(z) dz / (x_{k+1}^T - x_k^T)$ is the average value of $T(x)$ over the bay. The $T(x)$ is again obtained using spline curve fitting. Ideally, $u_{k+1} - u_k$ and $\int_{z=x_k^T}^{x_{k+1}^T} T(z) dz / (x_{k+1}^T - x_k^T)$ are equal. However, in this step they may never be equal, because the former is an integer while the latter can take any real values. Hence, the $u_{k+1}$ that minimizes their difference is selected.

After going through the length of the entire corridor, deselect the last transfer stop if it is too close to the first; i.e., if $\int_{z=x_{N^T}^T}^{L} \frac{dz}{s(z)} < \int_{z=x_{N^T}^T}^{L} T(z) dz / 2(L - x_{N^T}^T)$.

Step 3. Between any two consecutive transfer stops, the non-transfer stops are assigned to the $m_{cw}$ CW lines in a fixed order (as illustrated in Figure 1). Similarly assign non-transfer stops to the $m_{ccw}$ CCW lines. Note that the total number of non-transfer stops generated in steps 1-2 may not be divisible by $m_{cw}$ (or $m_{ccw}$). As a result, some lines may have one more non-transfer stops than other lines in the last skip-stop bay. However, the impact of this asymmetry on the generalized cost is trivial.

The three steps are summarized as follows:

**Algorithm 3: Generating a stop location plan.**

Use $s^*(x_j)$ and $T^*(x_j)$ to fit $s(x)$ and $T(x)$ to two spline functions, respectively.
Place one stop at every $x$ where $\int_{z=0}^{x} \frac{dz}{s(z)}$ is an integer and then obtain the stop location set $\Omega = \{x_i^S : i = 1, 2, \ldots, N^S\}$;
Remove the last stop from $\Omega$ and set $N^S \leftarrow N^S - 1$, if $L - x_{N^S}^S < \frac{s(L)}{2}$.
Select the first transfer stop at $x = 0$, i.e., $u_1 = 1$ and $x_1^T = 0$.
Set $k = 1$ and $N^T = 1$.
Do:
    For each $u_{k+1}$ in $\{u_k + 1, u_k + 2, \ldots, N^S\}$:
        If $u_{k+1} - u_k - 1$ is a common multiple of $m_{cw}$ and $m_{ccw}$:
            Let $x_{k+1}^T = x_{u_{k+1}}^S$;
            Record the optimal $u_{k+1}^*$ that minimizes $\left| (u_{k+1} - u_k) - \frac{\int_{z=x_k^T}^{x_{k+1}^T} T(z) dz}{x_{k+1}^T - x_k^T} \right|$.
        End If
    End For

$k \leftarrow k + 1, N^T \leftarrow N^T + 1$.
Until $u_k = N^S$
Deselect the $N^S$-th stop and set $N^T \leftarrow N^T - 1$.

Deselect the last transfer stop and set $N^T \leftarrow N^T - 1$ if $\int_{z=x_{N^T}^T}^{L} \frac{dz}{s(z)} < \frac{\int_{z=x_{N^T}^T}^{L} T(z) dz}{2(L - x_{N^T}^T)}$.

Between any two consecutive transfer stops, assign the non-transfer stops to the $m_{cw}$ CW lines and to the $m_{ccw}$ CCW lines in a fixed order.
**Output** the location and the type of each stop.

# Appendix F. Assessment of approximation errors

## F.1 An algorithm for calculating user and agency costs with exact stop locations

**Algorithm 4: Calculating user and agency costs using exact stop locations.**

Calculate the users' total access and egress time as $\int_{x=0}^{L}[P_{cw}(x) + Q_{cw}(x) + P_{ccw}(x) + Q_{ccw}(x)] \frac{g(x)}{v_w} dx$, where $g(x)$ denotes the distance between $x$ and the nearest stop.

For each $(i,j)$ in $\{1,2,\ldots,N^S\} \times \{1,2,\ldots,N^S\}$, where $N^S$ is the number of stops (both transfer and non-transfer ones):

Calculate for the trips from stop $i$ to stop $j$: (i) the aggregate demand, $\lambda_{i,j}$; (ii) direction of travel, $\rho_{i,j} \in \{CW, CCW\}$; (iii) the travel distance per trip, $d_{i,j}$; and (iv) type of the trips, $\gamma_{i,j} \in \{1,2,3,4,5\}$ (see Section 2.3.2 for the definition of trip types). Denote $n_{i,j}$ as the number of stops visited when traveling from stop $i$ to stop $j$ by the skip-stop service in direction $\rho_{i,j}$.

Case on $\gamma_{i,j}$:

Case "1"
Calculate $n_{i,j}$ as the number of stops visited along the trip by taking any skip-stop line in direction $\rho_{i,j}$. Calculate the total wait and in-vehicle travel time as $\left(\frac{H_{\rho_{i,j}}}{2} + \tau \cdot n_{i,j} + \frac{d_{i,j}}{v}\right) \lambda_{i,j}$.

Case "2" and "3"
Calculate $n_{i,j}$ as the number of stops visited along the trip by taking the line that visits both stops $i$ and $j$. Calculate the total wait and in-vehicle travel time as $\left(\frac{m_{\rho_{i,j}} H_{\rho_{i,j}}}{2} + \tau \cdot n_{i,j} + \frac{d_{i,j}}{v}\right) \lambda_{i,j}$.

Case "4"
Calculate $n_{i,j}$ as the number of stops visited along the backtracking trip. The trip contains a transfer that is made at the nearest transfer stop to both stops $i$ and $j$. Calculate the total wait, in-vehicle travel time, and transfer penalty as $\left(\frac{m_{cw}H_{cw} + m_{ccw}H_{ccw}}{2} + \tau \cdot n_{i,j} + \frac{d_{i,j}}{v} + C_t\right) \lambda_{i,j}$.

Case "5"
Calculate $n_{i,j}$ as the number of stops visited along the trip by taking first the line visiting stop $i$ and then transferring to the line visiting stop $j$. The transfer is assumed to be made at the first encountered transfer stop. Calculate the total wait, in-vehicle travel time, and transfer penalty as $\left(m_{\rho_{i,j}} H_{\rho_{i,j}} + \tau \cdot n_{i,j} + \frac{d_{i,j}}{v} + C_t\right) \lambda_{i,j}$.

End Case

Sum up the users' total access and egress time, wait, in-vehicle travel time, and transfer penalty to obtain the total user cost.

Calculate the agency cost as $\frac{\pi_k L}{\mu}\left(\frac{1}{H_{cw}} + \frac{1}{H_{ccw}}\right) + \frac{\pi_h}{\mu H_{cw}}\left(\tau(n_{cw} + 1) + \frac{L}{v}\right) + \frac{\pi_h}{\mu H_{ccw}}\left(\tau(n_{ccw} + 1) + \frac{L}{v}\right) + \frac{2\pi_i L}{\mu} + \frac{\pi_s}{\mu} N^S$, where $n_{cw}$ and $n_{ccw}$ denote the numbers of stops visited by any CW and CCW skip-stop line, respectively.

**Output** the sum of the user and agency costs.

## F.2 Errors in generalized cost and cost components

The generalized costs and cost components calculated using Algorithm 4 are compared against those of the heuristic solutions calculated by (14a) for the 144 numerical cases (see Section 4.1). The percentage errors are summarized in Table F1. They reveal how accurate our CA model is. Note that the generalized cost error never exceeds 1.2% and averages only 0.2%. Most errors between detailed cost components are very small too.

Table F1. Percentage cost errors between the heuristic solutions and the designs with exact stop locations

| Cost items | Average error (%) | Maximum error (%) |
|---|---|---|
| Generalized cost, $GC$ | 0.2% | 1.2% |
| User cost, $UT_a + UT_w + UT_v + UT_t$ | 0.3% | 1.4% |
| Agency cost, $AC_K + AC_H + AC_I + AC_S$ | 0.1% | 0.3% |
| Access/egress cost, $UT_a$ | 0.5% | 1.3% |
| Wait cost, $UT_w$ | 0.5% | 1.7% |
| In-vehicle travel cost, $UT_v$ | 0.4% | 2.0% |
| Transfer penalty, $UT_t$ | 1.2% | 6.8% |
| Distance-based operating cost, $AC_K$ | 0.0% | 0.0% |
| Time-based operating cost, $AC_H$ | 0.1% | 0.3% |
| Line infrastructure cost, $AC_I$ | 0.0% | 0.0% |
| Stop infrastructure cost, $AC_S$ | 0.5% | 1.2% |